\documentclass[12pt]{amsart}

\newtheorem{theorem}{Theorem}[section]
\newtheorem{lemma}[theorem]{Lemma}

\theoremstyle{definition}

\theoremstyle{remark}

\usepackage{amscd,amssymb}
\begin{document}

\baselineskip 13pt

\title[Test elements, retracts and automorphic orbits]
{Test elements, retracts and\\ automorphic orbits}
\author{Sheng-Jun Gong}
\address{Department of Mathematics, The University of Hong Kong, Pokfulam Road,
Hong Kong SAR, China} \email{sjgong@hkusua.hku.hk}
\thanks{Sheng-Jun Gong was partially supported by a University of Hong Kong Postgraduate Studentship}

\author{Jie-Tai Yu}\address{Department of Mathematics, The University of Hong Kong, Pokfulam Road,
Hong Kong SAR, China} \email{yujt@hkucc.hku.hk,    yujietai@yahoo.com}
\thanks{Jie-Tai Yu was partially supported by an RGC-CERG Grant}



\maketitle

\noindent {\bf Abstract.} Let $A_2$ be a
free associative or polynomial algebra of rank two over a field $K$ of characteristic zero.
 Based on the
degree estimate of Makar-Limanov and J.-T.Yu,
we prove: 1) An  element $p \in A_2$
is a test element if $p$ does not belong to any proper retract of $A_2$; 2) Every endomorphism preserving the
automorphic orbit of a nonconstant element of $A_2$ is an automorphism.

\baselineskip 13.5pt

\section{\large Introduction and main results}

\

In the sequel,  $K$  always denotes a field of of characteristic zero. Automorphisms (endomorphisms)
always mean $K$-automorphisms ($K$-endomorphisms).

\

Let $A_n$ be a free associative or polynomial algebra of rank $n$
over $K$. An element $p\in A_n$ is called a \emph{test element} if
every endomorphism of $A_n$ fixing $p$ is an automorphism. A
subalgebra $R$ of $A_n$ is called a \emph{retract} if there is an
idempotent endomorphism $\pi (\pi^2=\pi)$ of $A_n$ (called a
\emph{retraction} or a \emph{projection}) such that $\pi(A_n)=R$.
Test elements and retracts of groups and other algebras  are defined
in a similar way. Test elements and retracts of algebras and groups
have recently been studied in \cite{C, DY1, DY2, ES,  Je2, MSY2,
MUY2, MY1, MY2, MY3, MZ, Sh1, Sh2, SY1, SY2}.

\

A test element does not belong to any proper retract for any algebra or group as the corresponding non-injective idempotent endomorphism is not an automorphism.
The converse is proved by Turner \cite{T}  for free groups,
by Mikhalev and Zolotykh \cite{MZ} and by Mikhalev and J. -T. Yu \cite{MY1, MY2} for free Lie algebras and free Lie superalgebras respectively,
and by  Mikhalev, Umirbaev and J. -T. Yu \cite{MUY1}
for free nonassociative algebras.  See also Mikhalev, Shpilrain and J. -T. Yu \cite{MSY2}.

\

In view of the above, we may raise the
following

\

\noindent {\bf Conjecture 1.} If an element $p\in A_n$ does not belong to
any proper retract of $A_n$, then $p$ is a test element.

\

Recently, V.Shpilrain and J. -T. Yu \cite{SY2} proved
Conjecture 1  for $\mathbb{C}[x,y]$.
A key lemma in their proof is the degree estimate of Shestakov and Umirbaev \cite{SU1},
which plays a crucial role in the recent celebrated solution of the Nagata conjecture \cite{SU2, SU3} and the Strong Nagata conjecture \cite{UY}.

\

More recently, Makar-Limanov and J. -T. Yu \cite{MLY} developed a new combinatorial method based on the Lemma on radicals and obtained a sharp degree estimate for the `free' case, namely, for a free associative algebra or a polynomial algebra over a field of  characteristic
zero. It has found  applications for
automorphisms and coordinates of polynomial and free associative
algebras. See S.-J. Gong and J.-T. Yu \cite{GY2}.

\

Now we consider another related problem.
In an algebra or a group, certainly
an automorphism preserves the automorphic orbit of an element $p$. The converse is proved  by Shpilrain \cite{Sh3} and Ivanov \cite{I} for
free groups of rank two,  by D. Lee \cite{L} for free groups of any rank,
by Mikhalev and J. -T. Yu \cite{MY2} for free Lie algebras and by Mikhalev, Umirbaev and J.-T. Yu  \cite{MUY1} for free non-associative algebras, by van den Essen and Shpilrain \cite{ES}
for $A_2$ when $p$ is a coordinate, by Jelonek \cite{Je1} for polynomial
algebras over $\mathbb{C}$  when $p$ is a coordinate.
For the related linear coordinate preserving problem, see, for instance, S.-J. Gong and J.-T. Yu \cite{GY1}.
See also the  book \cite{MSY2}.

\

In view of the above, we may raise the
following

\

\noindent {\bf Conjecture 2.} Let $p\in A_n-K$.
Then any endomorphism of $A_n$ preserving
the automorphic orbit of $p$ must be an automorphism.

\

Conjecture 2 has recently been settled affirmatively by
J.-T. Yu \cite{Y} for $A_2=\mathbb{C}[x,y]$ based on Shpilrain and J.-T.Yu's characterization
of test elements of $\mathbb{C}[x,y]$ in \cite{SY2} and the main
result in Drensky and J.-T.Yu \cite{DY2}.

\

In this paper, based on the
recent degree estimate of Makar-Limanov and J.-T.Yu \cite{MLY},
the main ideals and techniques in
Drensky and J.-T.Yu \cite{DY2},
Shpilrain and J.-T. Yu \cite{SY1,SY2}, and
J.-T.Yu \cite{Y},
we prove both Conjecture 1 and Conjecture 2 for $n=2$. Our main results are

\

\begin{theorem}\label{test}  If an element $p\in A_2$ does not belong to any proper retract of $A_2$, then $p$ is a test element of $A_2$.
\end{theorem}

Theorem \ref{test} was proved by Shpilrain and J.-T.Yu \cite{SY2}
for $A_2=\mathbb{C}[x,y]$.

\begin{theorem}\label{orbit}
If an endomorphism $\phi$ of $A_2$ preserves
the automorphic orbit of a nonconstant element $p\in A_2$, then $\phi$ is an automorphism
of $A_2$.\end{theorem}

Theorem \ref{orbit} was proved by J.-T.Yu \cite{Y} for $A_2=\mathbb{C}[x,y]$.

\

Crucial to the proofs of the above two theorems are the following two results, which have their own interests.

\begin{theorem}\label{injection}
Let $p\in A_2$ have outer rank two.
Then any injective endomorphism $\phi$ of $A_2$
is an automorphism if $\phi(p)=p$.
\end{theorem}

Theorem \ref{injection} may be viewed
as an analogue of a result in Turner \cite{T} for free groups.
It was proved for $A_2=\mathbb{C}[x,y]$ in J.-T.Yu \cite{Y}
based on a result in Shpilrain and J.-T.Yu \cite{SY2}.

\begin{theorem}\label{fixed}
An element $p(x,y)\in A_2$ belongs to a proper retract of $A_2$  if $p(x,y)$ is fixed by a noninjective endomorphism $\phi$ of $A_2$. Moreover,
in this case there exists a positive integer $m$ such that $\phi^m$ is a retraction
of $A_2$.
\end{theorem}

Theorem \ref{fixed} was proved for $A_2=\mathbb{C}[x,y]$ in Drensky and J.-T.Yu \cite{DY2}.

\section{\large Proofs}

The following two lemmas are Theorem 1.1 and Proposition 1.2
in Makar-Limanov and J.-T.Yu \cite{MLY}.

\begin{lemma}\label{MLY1}
Let $A_n=K\langle x_1, \cdots, x_n \rangle$ be a free associative algebra over a field $K$ of characteristic zero, $f, g \in A$ be algebraically independent, $f^{+}$ and $g^{+}$ are algebraically independent, or $f^{+}$ and $g^{+}$ are algebraically dependent and neither $ \deg (f) \mid \deg (g)$ nor $ \deg (g) \mid \deg (f)$, $p \in K\langle x,y \rangle$. Then

$$\deg(p(f,g)) \ge \frac{\deg[f, g]}{\deg(fg)}w_{\deg(f), \deg(g)}(p).$$
\end{lemma}

Here $\deg$ is the total degree, $w_{\deg(f),\deg(g)}(p)$ is the weighted degree of $p$ when the weight of the first variable is $\deg(f)$ and the weight of the second variable is $\deg(g)$, $f^{+}$ and $g^{+}$ are the highest homogeneous components of $f$ and $g$ respectively, and $[f,g]=fg-gf$ is the commutator of $f$ and $g$.

\begin{lemma}\label{MLY2}
Let $A_n=K[x_1, \cdots, x_n]$ be a polynomial algebra over a field $K$ of  characteristic
zero, $f, g \in A$ be algebraically independent, $p \in K[x, y]$. Then

\

$\deg(p(f,g)) \ge$
$$w_{\deg(f), \deg(g)}(p)[1-\frac{(\deg(f), \deg(g))(\deg(fg)-\deg(J(f,g))-2)}{\deg(f)\deg(g)}].$$
\end{lemma}

Here $\deg$ is the total degree, $w_{\deg(f), \deg(g)}(p)$ is the weighted degree of $p$ when the weight of the first variable is $\deg(f)$ and the weight of the second variable is $\deg(g)$, $(\deg(f), \deg(g))$ is the greatest common divisor of $\deg(f)$ and $\deg(g)$, $\deg(J(f,g))$ is the largest degree of non-zero Jacobian determinants of $f$ and $g$ with respect to two of $x_1, \cdots, x_n$.

\

The following characterization of
a proper retract of $A_2$ was obtained by Shpilrain and J.-T.Yu \cite{SY1}
based on a result of Costa \cite{C}.

\begin{lemma}\label{retract}
Let $R$ be a proper retract of $A_2$. Then
$R=K[r]$ for some $r\in A_2$. Moreover, there exists an automorphism
$\alpha$ of $A_2$ such that $\alpha(r)=x+w(x,y)$,
where $w(x,y)$ belongs to the ideal of $A_2$ generated by $y$.
\end{lemma}

\begin{lemma}\label{orank} Let $p\in A_2$ with outer rank $2$ and $f, g\in A_n$.
Then $w_{\deg(f), \deg(g)}(p)\ge \deg(f)+\deg(g)$.
If every monomial of $p$ contains
both $x$ and $y$ and $\deg(p)>2$, then
$w_{\deg(f), \deg(g)}(p)> \deg(f)+\deg(g)$.
\end{lemma}

\begin{proof} 1) If $p$ contains a monomial containing both $x$ and $y$,
where $i\ne 0, j\ne 0$, $w_{\deg(f), \deg(g)}(p)\ge i(\deg(f))+j(\deg(g))\ge
\deg(f)+\deg(g)$. If every monomial of $p$ contains
both $x$ and $y$ and $\deg(p)>2$, then the second inequality
becomes strict.

\

2) Otherwise $p$ must contain monomials $x^i$ and $y^j$
where $i\ge 2,\ j\ge 2$. Then $w_{\deg(f), \deg(g)}(p)\ge 2\max\{\deg(f), \deg(g)\}
\ge\deg(f)+\deg(g)$.
\end{proof}

\begin{lemma}\label{commutator}
Let $A_n=K\langle x_1, \cdots, x_n \rangle$ be a free associative algebra over an arbitrary field $K$ of zero characteristic, $f, g \in A_2$ be algebraically independent,  $p \in K\langle x,y \rangle$ have outer rank two. Then
$$\deg(p(f,g)) \ge \deg[f,g].$$
If every monomial of $p$ contains
both $x$ and $y$ and $\deg(p)>2$, then $$\deg(p(f,g))>\deg[f,g].$$
\end{lemma}
\begin{proof} Let
1) If $f^+$ and $g^+$ are algebraically independent; or $f^+,\ g^+$
are algebraically dependent, but
$\deg(f)\nmid \deg(g)$ and $\deg(g)\nmid \deg(f)$. Then by
Lemma \ref{MLY1} and Lemma \ref{orank},
$\deg(p(f,g)) \ge \deg[f,g]$. If, in addition,
every monomial of $p$ contains
both $x$ and $y$ and $\deg(p)>2$, then by Lemma \ref{MLY1} and Lemma \ref{orank},
$\deg(p(f,g))>\deg [f,g]$.

\

2) Otherwise there exists an automorphism $\alpha$,
which is the composition of a sequence of
elementary automorphisms, such that $\alpha(f)=\bar f,\
\alpha(g)=\bar g$,\ $\bar p =\alpha^{-1}(p)$ satisfying
the condition in 1). Then
$\deg(p(f,g))$

\noindent $=\deg(\bar p(\bar f, \bar g))\ge\deg[\bar f, \bar g]
=\deg [f,g]$.
\end{proof}

\begin{lemma}\label{jacobian}
Let $A_n=K[x_1, \cdots, x_n]$ be a polynomial algebra over an arbitrary field $K$ of zero characteristic, $f, g \in A_n$ be algebraically independent, $p \in K[x, y]$ has outer rank two. Then
$$\deg(p(f,g)) \ge\deg(J(f,g))+2.$$
\end{lemma}
\begin{proof} We may assume $\deg(f)=m, \deg(g)=n$.
As $p$ has outer rank $2$, by Lemma \ref{orank}
then  $p$ contains a monomial with both $x$ and $y$,
or contains monomials $x^i$ and $y^j$ where
$i\ge 2,\ j\ge 2$.

\

1) Let $f^+$ and $g^+$ be algebraically
independent.

   a) If there exists a monomial in $p$
   containing both $x$ and $y$, then
   $\deg(p(f,g))\ge\deg(f)+\deg(g)\ge\deg(J(f,g))+2$;

   b) Otherwise $p$ must have a monomial of $x^i$ where $i\ge 2$,
   and another monomial $y^j$
      where $j\ge 2$, then $\deg(p(f,g))\ge 2\max\{m,n\}\ge\deg(f)+\deg(g)\ge\deg(J(f,g))+2$;

\

2) Let $f^+,\ g^+$ be
 algebraically dependent, and
$m\nmid n$ and $n\nmid m$.

   c) If $w_{\deg(f),\deg(g)}(p)<$ lcm$(m,n)$,
   then in $p(f,g)$, $f^+$ and $g^+$
   cannot cancel out, hence similar to the case 1 a),  $\deg(p(f,g))\ge\deg(f)+\deg(g)
   \ge\deg(J(f,g))+2$.

   d) Otherwise $w_{\deg(f),\deg(g)}(p)\ge$lcm$(m,n)=mn/(m,n)$.
We also have $mn=(m,n)\text{lcm}(m,n)\ge(m,n)(m+n)$.
Hence $\deg(p(f,g))\ge\deg(J(f,g))+2$ by Lemma \ref{MLY2}.

\

3) Let $f^+,\ g^+$
be algebraically dependent, but
$m\mid n$ or $n\mid m$.
Then by same process in the Proof 2) of
Lemma 2.4, we may reduce to the above case 1) or case 2).
\end{proof}

\begin{lemma}\label{k}
Let $\phi=(f,g)$ be an injective endomorphism of
$K\langle x,y\rangle$ but not
an automorphism. Then
$\deg([\phi^k(x), \phi^k(y)])\ge k+2$ for $k\ge 0$.
\end{lemma}
\begin{proof} We use induction.
$\deg[\phi^0(x), \phi^0(y)]=\deg[x,y]=0+2$.
Assuming $\deg[\phi^{k-1}(x),\phi^{k-1}(y)]\ge (k-1)+2.$
Define $p(x,y):=[f(x,y), g(x,y)]$.
As $\phi=(f,g)$ is injective, every monomial of $p(x,y)$ contains both $x$ and $y$.
Since
$\phi=(f,g)$ is not an automorphism,
by the well-known result
of Dicks (see, Dicks \cite{D}, or Cohn \cite{Cn}),
$\deg(p(x,y))>\deg(x)+\deg(y)=2$.
Applying Lemma \ref{commutator},
$\deg(p(u,v))>\deg[u,v]$\  for $u=\phi^{k-1}(x),\ v=\phi^{k-1}(y)$,\
hence

\noindent $\deg[\phi^k(x),\phi^k(y)]=\deg(p(\phi^{k-1}(x),\phi^{k-1}(y)))
>\deg[\phi^{k-1}(x),\phi^{k-1}(y)]\ge (k-1)+2=k+1$.
Therefore, $\deg[\phi^k(x),\phi^k(y)]\ge(k+1)+1=k+2$.
\end{proof}

\begin{lemma}\label{k'}
Let $\phi=(f,g)$ be an injective endomorphism of
$K[x,y]$ but not an automorphism and there exists an element $p\in K[x,y]$
fixed by $\phi$. Then
$\deg(J(\phi^k(x),\phi^k(y)))\ge k$ for $k\ge 0$.
\end{lemma}
\begin{proof}
As $\phi$ fixes $p$,  $\phi$ is not an automorphism,
by a result of Kraft \cite{K} (see also Shpilrain and J.-T.Yu
\cite{SY1}), $\deg(J(\phi(x), \phi(y)))=\deg(J(f,g))$

\noindent $\ge 1$.
By the chain rule for the Jacobian,

\

\noindent $\deg(J(\phi^k(x),\phi^k(y)))$

\noindent $=\deg(J(f,g)(\phi^{k-1}(x),\phi^{k-1}(y))(J(\phi^{k-1}(x),\phi^{k-1}(y))))$

\noindent $\ge \deg(J(\phi^{k-1}(x),\phi^{k-1}(y)))+1.$

\

\noindent The proof is concluded by induction.

\end{proof}

\begin{lemma}\label{M}
Let $\phi=(f,g)$ be an injective endomorphism of
$A_2$ but not an automorphism.  Then any element $p\in A_2$
with outer rank $2$ cannot be fixed by $\phi$.
\end{lemma}
\begin{proof} If $p\in A_2$ with outer rank two fixed by $\phi$, then
$\deg(p(f,g))=\deg(p(\phi^k(x), \phi^k(y)))\ge k+2$ for all $k\ge 0$.
by Lemma \ref{commutator} and Lemma \ref{k} for noncommutative case;
and by Lemma \ref{jacobian} and Lemma \ref{k'} for polynomial case.
The contradiction completes the proof.
\end{proof}

\noindent {\bf Proof of Theorem \ref{injection}.}

\

By Lemma \ref{M}.$\Box$

\

\noindent {\bf Proof of Theorem \ref{fixed}.}

\

The proof presented here is similar to
the proof of the main Theorem in Drensky and J.-T. Yu \cite{DY2}.

Let $p\in A_2-\{0\}$
fixed by a noninjective endomorphism of $A_2$. Then
$\phi(x)$ and $\phi(y)$ are algebraically dependent over $K$.
Let us denote the image of $\phi(A_2)$ by $S=K[\phi(x), \phi(y)]$
 (since $\phi(x)$ and $\phi(y)$ are algebraically dependent,
 $\phi(x)$ and $\phi(y)$ are in a polynomial algebra of rank one over $K$
 as a consequence of a result of Bergman \cite{B} for noncommutative case
 and as a consequence of a result of Shestakov and Umirbaev \cite{SU1}
 for polynomial case) and
by $Q(S)$ the field of fractions of $S$.  Therefore the transcendence degree of $Q(S)$ over $K$
is 1. Let $0\ne q(x,y)\in(\text{\rm Ker}(\phi))\cap S$. Since $p(x,y)$
also belongs to $S$, the polynomials $p$ and $q$ are algebraically dependent and
$$
h(p,q)=a_0(q)p^n+a_1(q)p^{n-1}+\ldots+a_{n-1}(q)p+a_n(q)=0
$$
for an irreducible polynomial $h(u,v)\in K[u,v]$ and
$a_i(t)\in K[t]$, $i=0,1,\ldots,n$. Hence
$\phi(h(p,q))=h(\phi(p),\phi(q))=h(p,0),$
$$a_0(0)p^n+a_1(0)p^{n-1}+\ldots+a_{n-1}(0)p+a_n(0)=0.$$
Therefore $a_0(0)=a_1(0)=\ldots=a_n(0)=0$. Now the polynomials
$a_i(t)$ have no constant terms and $h(u,v)$ is divisible by $v$ which
contradicts to the irreducibility of $h(u,v)$. Therefore
$(\text{\rm Ker}(\phi))\cap S=0$ and $\phi$ acts injectively on its image $S$.
Hence we may extend the action of $\phi$ on $Q(S)$ (because
$a_1/b_1=a_2/b_2$ in $Q(S)$ is equivalent to $a_1b_2=a_2b_1$ and hence
$\phi(a_1/b_1)=\phi(a_1)/\phi(b_1)=\phi(a_2)/\phi(b_2)=\phi(a_2/b_2)$).
By L\"uroth's theorem  (See, for instance, Schinzel \cite{Sc}),
$Q(S)=K(w)$ for some $w\in Q(S)$.
The automorphism $\phi$ fixes $p(x,y)$ and its extension $\bar\phi$ on $Q(S)$ fixes
$K(p)$. Since $w$ is algebraic over $K(p)$,\
$Q(S)$ is a finite dimensional vector space over $K(p)$ and
$\bar\phi$ is a $K(p)$-linear operator of $Q(S)$ with trivial kernel.
Hence $\bar\phi$ is invertible on $Q(S)$ and we may consider $\bar\phi$ as
an automorphism of the finite field extension $Q(S)$ over $K(p)$ which fixes
$K(p)$. By Galois theory ($\bar\phi$ interchanges the
roots of the minimal polynomial of $w$ over $K(p)$
and there are finite number of possibilities for $\bar\phi(w)$),
$\bar\phi$ has finite order. Let $\bar\phi^m=1$. Then
$\phi^{m+1}(r)=\phi^m(\phi(r))=\bar\phi^m(\phi(r))=\phi(r)$
for every $r\in A_2$ and
$(\phi^m)^2=\phi^{m+1}\phi^{m-1}=\phi\phi^{m-1}=\phi^m$. Therefore
$\pi=\phi^m$ is a retraction (idempotent endomorphism) of $A_2$ with a nontrivial kernel and
$\pi(p)=p$. Hence $p(x,y)$ is in the image of $\pi$ which is a proper retract
$\pi(A_2)$ of $A_2$. $\Box$

\

\noindent {\bf Proof of Theorem \ref{test}.}

\

 As $p\in A_2$ does not belong to
 any proper retract of $A_2$, by Theorem \ref{fixed}, any
 endomorphism $\phi$ of $A_2$ fixing
 $p$ must be injective. By Lemma \ref{retract}, obviously $p$ must have outer rank two,
 otherwise $p$ would belong to a proper retract of $A_2$. By Theorem \ref{injection},
 $\phi$ is an automorphism. Hence $p$ is a test element
 of $A_2$. $\Box$

 \

\noindent {\bf Proof of Theorem \ref{orbit}.}

\

The proof presented here is similar to
the proof of the main result Theorem 1.4 in J.-T. Yu \cite{Y}.

\

We may assume that $\phi(p)=p$. By the definition
of the test element, we may assume $p$ is not a test element. By Theorem \ref{test}, we may assume $p$ belongs to a proper retract $K[r]$ of $A_2$. By a result in J.-T.Yu \cite{Y}, we may assume $p$ has outer rank $2$. By Theorem \ref{injection}, we may assume $\phi$ is non-injective. Suppose that $p=f(r)$, where $f \in K[t]-K$, $\deg(f)=m$. By Theorem \ref{fixed}, $\pi=\phi^m$ is a retraction of $A_2$ to $K[r]$. As $\phi$ preserves the automorphic orbit of $p$, so does $\pi=\phi^m$. Applying Lemma \ref{retract}
(suppose $\alpha(r)=x+w(x,y)$, where $w(x,y) \notin K[y]$ belongs to the ideal of $A_2$ generated by $y$, $\alpha$ is some automorphism of $A_2$, replace $r$ by $\alpha(r)$, and $\pi$ by $\alpha \pi \alpha^{-1}$), we have reduced our proof to the following

\

\begin{lemma}
Let $r=x+w(x,y)$, where $w(x,y)$ belongs to the ideal of $A_2$ generated by $y$ and $w(x,y) \notin K[y]$, $\pi$ the retraction of $A_2$ onto $K[r]$ defined by $\pi(x)=x+w(x,y)$, $\pi(y)=0$, $f \in K[t]-K$. Then $\pi$ does not preserve the automorphic orbit of $f(r)$.
\end{lemma}

\begin{proof}

Suppose on the contrary, $\pi$ preserves the automorphic orbit of $f(r)$. Then for any automorphism $\alpha$ of $A_2$, $\pi\alpha(f(r))=\beta(f(r))\in K[r]$ for some automorphism $\beta$ of $A_2$. Note that $\pi\beta(f(r))=\pi^2\beta(f(r))=\pi\alpha(f(r))=\beta(f(r))$. By Theorem \ref{fixed}, $\pi^{\deg(f)}=\pi$ is the retraction of $A_2$ onto the retract $K[\beta(r)]$ taking $\beta(r)$ to $\beta(r)$. By hypothesis, $\pi$ is also a retraction of $A_2$ onto the retract $K[r]$ taking $r$ to $r$. This forces that $\beta(r)=cr+d$ for some $c \in K^*,\ d\in K$. We have concluded that for any automorphism $\alpha$ of $A_2$, there exists some $c \in K^*,\ d\in K$, such that $\pi \alpha(f(r))=f(cr+d)$.

\

Now we proceed the proof in two cases.

\

\noindent {\bf 1. Noncommutative case: $A_2=K\langle x,y\rangle$.}

\

Denote by $\mathcal{C}$ the commutator ideal of $K\langle x, y \rangle$.

\

 a) If $w(x,y) \in \mathcal{C}$, then take $\alpha$ to be the automorphism of $K\langle x, y \rangle$ defined by $\alpha(x)=y+x^2$, $\alpha(y)=x$. Direct calculation shows that $\pi \alpha(f(r))=f(r^2+w(r^2, r))=f(r^2) \neq f(cr+d)$, a contradiction.

\

b) If $w(x,y) \notin \mathcal{C}$, then $w^{a}(x,y)=yv(x, y)$ for some $v(x,y) \in K[x,y]-\{0\}$. Here $w^{a}(x,y) \in K[x,y]$ is the image of $w(x,y)$ under the  abelianization from  $K\langle x, y \rangle$ onto $K[x,y]$. Let $M$ be a positive integer greater than $\deg (v(x,y))$, it is easy to see that $x^M-y$ does not divide $v(x,y)$ in $K[x, y]$. Let $\alpha$ be the automorphism of $K\langle x, y \rangle$ defined by $\alpha(x)=x$, $\alpha(y)=y+x^M$. Then  $\pi \alpha(f(r))=f(r+w(r, r^M))=f(r+r^{M}v(r, r^M))$. As $x^M-y$ does not divide $v(x,y)$, $v(r, r^M) \neq 0$. Therefore $\pi \alpha(f(r))=f(r+r^{M}v(r, r^M)) \neq f(cr+d)$, a contradiction.

\

\noindent {\bf 2. Polynomial case: $A_2=K[x,y]$.}

\

In this case we write $w(x,y)=yq(x,y)$ where $q(x,y)\notin K[y]$.
Let $M$ be a positive integer greater than $\deg(q(x,y))$,
it is easy to see that $x^M-y$ does not divide $q(x,y)$ in
$K[x,y]$. Let $\alpha$ be the automorphism
of $K[x,y]$ defined by $\alpha(x)=x,~ \alpha(y)=y+x^M$.
Then easy calculation shows that $\pi\alpha(f(r))=f(r+r^Mq(r,r^M))$.
As  $x^M-y$ does not divide $q(x,y)$,\ $q(r,r^M)\ne 0$.
Therefore $\pi\alpha(f(r))=f(r+r^Mq(r,r^M))\ne f(cr+d)$. The contradiction completes
the proof. \end{proof}

\

\section{\large Acknowledgements}

The authors are grateful to the referee for very helpful comments,
in particular, for improvement of the proof of Lemma \ref{k}.

\

\end{document}